\documentclass[%
reprint,
nofootinbib,
amsmath,amssymb,
aps,
showkeys,
]{revtex4-1}

\usepackage{graphicx}% Include figure files
\usepackage{dcolumn}% Align table columns on decimal point
\usepackage{bm}% bold math
\usepackage{ dsfont }
\usepackage{nicefrac}
\usepackage{algorithm}
\usepackage{algpseudocode}
\usepackage{accents}
\newcommand{\ubar}[1]{\underaccent{\bar}{#1}}

\begin{document}
	
	\preprint{APS/123-SDE}
	
	\title{Probabilistic simulation of partial differential equations}% Force line breaks with \\DFI / 
	
	\author{Philipp Frank}
	%\email{philipp@mpa-garching.mpg.de}
	\author{Torsten A. En\ss lin}%
	\affiliation{%
		Max-Planck Institut f{\"u}r Astrophysik, Karl-Schwarzschild-Str. 1, 85748, Garching, Germany
	}%
	\affiliation{%
		Ludwig-Maximilians-Universit{\"a}t M{\"u}nchen, Geschwister-Scholl-Platz 1, 80539, M{\"u}nchen, Germany
	}%

	\date{\today}% It is always \today, today,

\begin{abstract}
Computer simulations of differential equations require a time discretization, which inhibits to identify the exact solution with certainty. Probabilistic simulations take this into account via uncertainty quantification. The construction of a probabilistic simulation scheme can be regarded as Bayesian filtering by means of probabilistic numerics. Gaussian prior based filters, specifically Gauss-Markov priors, have successfully been applied to simulation of ordinary differential equations (ODEs) and give rise to filtering problems that can be solved efficiently. This work extends this approach to partial differential equations (PDEs) subject to periodic boundary conditions and utilizes continuous Gaussian processes in space and time to arrive at a Bayesian filtering problem structurally similar to the ODE setting. The usage of a process that is Markov in time and statistically homogeneous in space leads to a probabilistic spectral simulation method that allows for an efficient realization. Furthermore, the Bayesian perspective allows the incorporation of methods developed within the context of information field theory such as the estimation of the power spectrum associated with the prior distribution, to be jointly estimated along with the solution of the PDE.
\end{abstract}

%Use showkeys class option if keyword
\maketitle

%\tableofcontents

\section{Introduction}
Numerical simulation of partial differential equations (PDEs) has been studied extensively for a long time as PDEs arise naturally in many scientific fields.
Recently, fully probabilistic approaches to simulation have been proposed \cite{RAISSI2017683, doi:10.1137/17M1120762}, many of them within the context of probabilistic numerics (PN) \cite{cockayne2017probabilistic, KerstingMahsereci2020}. Many probabilistic numerical methods aim to disentangle traditional numerical algorithms into the prior assumptions as well as the (artificial) observations that appear within the algorithm \cite{SchoberDH2014}.
This provides an uncertainty quantification within the context of Bayesian reasoning \cite{doi:10.1098/rspa.2015.0142} and often has led to new variants of the algorithms by replacing prior assumptions \cite{KerHen16}.

In this work we aim to discuss probabilistic numerical simulation within the context of information field theory (IFT) \cite{doi:10.1063/1.4819999}, that is information theory for quantities that are defined over continuous spaces (i.E. fields). Previous works towards an information field theoretical consideration of PDE simulation has been established by means of information field dynamics (IFD) \cite{PhysRevE.87.013308, PhysRevE.97.033314}. IFD aims to construct a simulation step that is optimal in the information theoretical sense, that is minimal loss of information about the system between subsequent simulation steps. In this work, however, we follow a line of argument more closely related to PN rather then IFD. We discuss the relations to IFD in further detail once we established the main properties of the probabilistic solver.
Nevertheless, the usage of IFT allows for an application of non-parametric estimation of power spectra \cite{PhysRevD.83.105014} to the task of PDE simulation. This enables us to construct more sophisticated filters that adapt to the correlation structure of the simulated process.

We notice that our approach has considerable structural similarities to a recent reformulation of probabilistic simulation of ordinary differential equations (ODEs) by means of nonlinear Bayesian filtering \cite{Tronarp2019}, however here applied to PDEs.

\subsection{Introduction to IFT and notation}
In IFT we consider fields $s^x$ that are defined over some continuous domain $\Omega \subset \mathds{R}^d$ where $d$ denotes the dimensionality of the space and $x$ may label a location in a coordinate system on $\Omega$. We aim to provide probabilistic reasoning for fields, and therefore we need to define probability distributions for fields.
To this end we equip the function space $L^2\lbrace\Omega\rbrace$ with a scalar product defined as
\begin{equation}
a^\dagger b \equiv \int_{\Omega} a^{*}_x \ b^x \ \mathrm{d}x \ ,
\end{equation}
where $*$ denotes complex conjugation.
Consequently, applications of linear operators $O : L^2\lbrace\Omega\rbrace \rightarrow L^2\lbrace\Omega\rbrace$ are denoted as
\begin{equation}\label{eq:lincont}
b^x = (O a)^x = O^x_{\; x'} \ a^{x'} = \int_{\Omega} O^x_{\; x'} \  a^{x'} \ \mathrm{d}x' \ ,
\end{equation}
where we also introduced the continuous version of the Einstein sum convention.
This allows us to define a Gaussian distribution with mean $m$ and covariance $D$ for a field $s$ via
\begin{align}\label{eq:fieldgauss}
P(s) &= \mathcal{G}(s-m,D) \notag\\ &\equiv \frac{1}{\left|2 \pi D \right|^{\frac{1}{2}}} \ e^{-\frac{1}{2} (s-m)^\dagger D^{-1} (s-m)} \ ,
\end{align}
where $| \bullet |$ denotes the functional determinant. (For further details see e.g. \cite{2019AnP...53100127E}).
In order to perform inference we additionally need to define a mapping $R : L^2\lbrace\Omega\rbrace \rightarrow \mathds{R}^N$ (often referred to as response, or design-matrix) that maps a field $s$ to some discrete measurement data $d \in \mathds{R}^N$. Similar to Eq.\ \eqref{eq:lincont} we write
\begin{equation}
	d^i = (R s)^i = R^i_{\; x} s^x = \int_{\Omega} R^i_{\; x} \  a^{x} \ \mathrm{d}x \ .
\end{equation}
If we aim to apply the adjoint of $R$ (denoted as $R^\dagger$), however, we get that
\begin{equation}
	b^x = (R^\dagger d)^x = (R^\dagger)^x_{\; i} d^i \equiv \sum_{i = 1}^{N} R^x_{\; i} \ d^i \ ,
\end{equation}
as we define the scalar product in discrete space as a sum.

\section{Probabilistic Simulation within IFT}
To summarize some key results of probabilistic simulation required for PDE simulation, we start with a brief discussion of ODE simulation and show its relation to Bayesian filtering. For an extensive overview of PN methods for ODE simulation please refer to \cite{SchoberSarkkaHennig2018, Tronarp2019}.
\subsection{Probabilistic ODE simulation}
Consider an ODE of the form
\begin{equation}\label{eq:ode}
	\dot{s}^t \equiv \frac{\partial s^t}{\partial t} = f\left(s^t\right) \quad \text{with initial condition} \quad s^{t_0} = s^0 \ ,
\end{equation}
where $s^t \in \mathds{R}^M$ denotes the state of the system at time $t$ and $f : \mathds{R}^M \rightarrow \mathds{R}^M$ is a (non-linear) map.

A Bayesian approach to simulation can the be formulated as:
Given some prior knowledge on the field $s$ given as $P(s|s^0)$ we aim to constrain this prior via artificial observations such that it solves Eq.\ \eqref{eq:ode}. The resulting posterior distribution is thus informed via the information in the observations, as well as the prior assumptions.
To realize the ODE constraints, we may define a continuous data-set $d^t$ as
\begin{equation}
	d^t = \dot{s}^t - f\left(s^t\right) \ ,
\end{equation}
and require that $d^t = 0 \quad \forall t$.
In general, however, this gives rise to an infinite set of non-tractable constraints and given only finite computational resources, leads to non-computable posterior distributions.
Therefore, in the spirit of PN, we require this constraint to be satisfied only at a discrete set of moments in time $T \equiv \left\lbrace t_i\right\rbrace_{i \in \left\lbrace 0, ..., N-1\right\rbrace}$ via
\begin{equation}
	d = R \left(\dot{s} - f\left(s\right)\right) \quad \text{with} \quad R^i_{\; t} = \delta\left(t_i - t\right) \ ,
\end{equation}
and then require $d^i = 0 \quad \forall i \in \left\lbrace 0, ..., N-1\right\rbrace$.
Note that the choice of $R$ has an impact on the resulting simulation scheme as it introduces a measure and consequently a PN method for simulation is only fully specified given a prior distribution of the continuous process, as well as a measurement operation. The specific choice of $R$ considered in this work has the desirable property that
\begin{equation}
	R f\left(s\right) = f\left(R s\right) \ .
\end{equation}
As it will turn put, this property allows us to set up a simulation scheme that only requires to construct the distribution of $R s$ and $R \dot{s}$ from the prior.

To do so, consider the special case of a Gaussian prior for $s$ of the form of Eq.\ \eqref{eq:fieldgauss}. Furthermore let
\begin{equation}
	x = \left(\begin{matrix}
	\bar{s} \\
	\dot{\bar{s}}
	\end{matrix}\right) \equiv \left(\begin{matrix}
	R s \\
	R \dot{s}
	\end{matrix}\right) = \left(\begin{matrix}
	R s \\
	R \partial_t s
	\end{matrix}\right) \ ,
\end{equation}
where $\partial_t$ denotes the derivative w.r.t.\ $t$.
As Gaussian distributions are closed under affine transformations, we get that $x$ is also Gaussian distributed with mean
\begin{equation}
	m_x = \left(\begin{matrix}
	R m \\
	R \dot{m}
	\end{matrix}\right) \ ,
\end{equation}
and covariance $X$
\begin{equation}
	X = \left(\begin{matrix}
	R D R^\dagger & R D \partial_t^\dagger R^\dagger \\
	R \partial_t D R^\dagger & R \partial_t D \partial_t^\dagger R^\dagger
	\end{matrix}\right) \ ,
\end{equation}
where  $\partial_t^\dagger$ denotes taking the derivative to the left (i.E. the second index of $D$ in this case).
We can use these results to construct the posterior distribution of $s$ given $d=0$. Let $\ubar{T} \equiv [t_0, \infty) \setminus T$ and let $\ubar{s}$ be all $s^t$ with $t \in \ubar{T}$, we get that
\begin{align}
	&P\left(s|d = 0, s^0\right) = \int \mathrm{d}\dot{\bar{s}} \ P\left(\ubar{s}, \bar{s}, \dot{\bar{s}} | d = 0\right) \\ &\propto \int \mathrm{d}\dot{\bar{s}} \ P\left(d = 0 | \ubar{s}, \bar{s}, \dot{\bar{s}}\right) \ P\left(\ubar{s} | \bar{s}, \dot{\bar{s}} \right) \ P\left(\bar{s}, \dot{\bar{s}} | s^0\right) \notag\\ &=
	\int \mathrm{d}\dot{\bar{s}} \ \delta\left(\dot{\bar{s}} - f\left(\bar{s}\right)\right) \ P\left(\ubar{s} | \bar{s}, \dot{\bar{s}} \right) \ P\left(\bar{s}, \dot{\bar{s}} | s^0\right)
	\notag\\ &= P\left(\ubar{s} \left| x = \left(\begin{matrix}
	\bar{s} \\
	\dot{\bar{s}} = f\left(\bar{s}\right)
	\end{matrix}\right) \right)\right. \left.P\left(x = \left(\begin{matrix}
	\bar{s} \\
	\dot{\bar{s}} = f\left(\bar{s}\right)
	\end{matrix}\right)\right| s^0\right) \ .
\end{align}
First, we notice that the posterior for all $\ubar{s}$ remains a Gaussian distribution irrespective of $f$ and is equal to the conditional distribution of $s$ given the values and the first derivatives at all $T$. Furthermore we may write
\begin{equation}\label{eq:postgen}
	\left.P\left(x = \left(\begin{matrix}
	\bar{s} \\
	\dot{\bar{s}} = f\left(\bar{s}\right)
	\end{matrix}\right)\right| s^0\right) = P\left(\dot{\bar{s}} = f\left(\bar{s}\right) | \bar{s}, s^0 \right) P\left(\bar{s}| s^0\right) \ ,
\end{equation}
which ultimately renders the task of simulation a non-linear Bayesian regression problem in $\bar{s}$ \cite{Tronarp2019}.

\subsubsection{Gauss-Markov processes}
For general Gaussian priors, i.E. for general $D$ (see Eq.\ \eqref{eq:fieldgauss}), this approach scales with $N^2$ ($N^3$ in case of unknown hyper-parameters in $D$) as we need to compute conditional distributions for all $T$.
Therefore, as proposed by e.g.\ \cite{Schober2018}, one can achieve linear scaling in $N$ via usage of Gauss-Markov processes. In this work we restrict ourselves to the simple case of an integrated Wiener process (IWP), however a generalization to higher order Gauss-Markov process priors is possible as provided by \cite{Schober2018}. The IWP may be defined as
\begin{equation}
	\ddot{s}^t = \sigma \xi^t \quad \text{with} \quad \xi \sim \mathcal{G}\left(\xi, \mathds{1}\right) \ ,
\end{equation}
and yields the conditional distribution for $s^t$ and $\dot{s}^t$ given their values at a previous time step:
\begin{align}\label{eq:iwpprior}
	&\left.P\left(\left(\begin{matrix}
	s^{t_i} \\ \dot{s}^{t_i}
	\end{matrix}\right) \right| \left(\begin{matrix}
	s^{t_{i-1}} \\ \dot{s}^{t_{i-1}}
	\end{matrix}\right)\right) \notag\\ &= \mathcal{G}\left(\left(\begin{matrix}
	s^{t_i} \\ \dot{s}^{t_i}
	\end{matrix}\right) - \left(\begin{matrix}
	1 & \Delta_i \\ 0 & 1
	\end{matrix}\right) \left(\begin{matrix}
	s^{t_{i-1}} \\ \dot{s}^{t_{i-1}}
	\end{matrix}\right) , \sigma^2 \left(\begin{matrix}
	\nicefrac{\Delta_i^3}{3} & \nicefrac{\Delta_i^2}{2} \\ \nicefrac{\Delta_i^2}{2} & \Delta_i
	\end{matrix}\right)\right) \ ,
\end{align}
where $\Delta_i = t_i - t_{i-1}$.

Using the IWP prior, the posterior Eq.\ \eqref{eq:postgen} reads
\begin{align}
	&P\left(\bar{s} | d, s^0\right) \propto \prod_{i=1}^{N-1} \left.P\left(\left(\begin{matrix}
	s^{t_i} \\ f\left(s^{t_i}\right)
	\end{matrix}\right) \right| \left(\begin{matrix}
	s^{t_{i-1}} \\ f\left(s^{t_{i-1}}\right)
	\end{matrix}\right)\right)  \notag\\ &= \prod_{i=1}^{N-1} \left[P\left(\dot{s}^{t_i} = f\left(s^{t_i}\right) | s^{t_{i}}, s^{t_{i-1}}, \dot{s}^{t_{i-1}} = f\left(s^{t_{i-1}}\right) \right)\right. \notag\\& \quad \left.P\left(s^{t_{i}} | s^{t_{i-1}}, \dot{s}^{t_{i-1}} = f\left(s^{t_{i-1}}\right) \right)\right] \ .
\end{align}
In words, the observations constructed via $R$ only affect the posterior locally and therefore the Markov property of the prior remains present in the posterior. As a consequence the Bayesian filtering problem defined in Eq.\ \eqref{eq:postgen} decomposes into a set of $N-1$ subsequent filtering problems, one for each $s^{t_i}$.

\subsection{PDEs with periodic boundary conditions}
To construct a probabilistic method for PDEs consider a generic PDE in $1 + 1$ dimensions for a scalar field $s$ of the form
\begin{equation}\label{eq:pde}
	\dot{s}^{t x} = f\left(s^{t x}, \left(s^{(1)}\right)^{t x}, \left(s^{(2)}\right)^{t x}, ... \right) \ , 
\end{equation}
with $f : \mathds{R} \otimes \mathds{R} \otimes ... \rightarrow \mathds{R}$, and $s^{(c)}$ denotes the $c$th spatial derivative of $s$. We restrict the discussion to scalar fields in $1 + 1$ dimensions but note that an extension to higher dimensions and vector fields is possible.
Furthermore we only consider PDEs that are compatible with periodic boundary conditions in the spatial domain\footnote{Other boundary conditions can be enforced by modification of the dynamical equations in the here presented approach, and possibly by a zero  padding area between those in the periodic domain. We leave this to future research.}, and, without loss of generality, require the size of the spatial domain to be equal to one.

For a probabilistic solver, we require a prior distribution for $s$. We remain in the setting of a Gauss-Markov prior and additionally assume independence of space and time prior correlations. I.e.\ we assume that
\begin{equation}
	\left< s^{t x} s^{t' x'}\right>_{P(s)} = C^{t t'} \ S^{x x'} = C^{t t'} S\left(|x - x'|\right) \ ,
\end{equation}
where we include the additional assumption that the spatial correlation structure is a priori statistical homogeneous and isotropic. We set $C$ such that $s$ follows an IWP in time.
Furthermore, we define $s$ in terms of its Fourier series
\begin{equation}
	s^{t x} = \sum_{k = -\infty}^{\infty} \tilde{s}^{t k} e^{2 \pi i k x} \ ,
\end{equation}
and use the fact that the Fourier modes $\tilde{s}$ of a statistically homogeneous process become statistically independent in Fourier space. The prior assumptions additionally imply that the time evolution of each Fourier mode $\tilde{s}^k$ follows an IWP of the form
\begin{align}\label{eq:contiwp}
	\ddot{\tilde{s}}^{t k} &= \sigma^{k} \ \xi^{t k} \quad \text{with} \quad \xi \sim \mathcal{G}\left(\xi,\mathds{1}\right) \ ,
\end{align}
with $\sigma$ such that $|\sigma|^2$ equals the Fourier spectrum associated with the spatial covariance $S$.

\subsubsection{Discrete Measurements}
In analogy to the ODE discussion we have to define a discrete set of measurements in order to arrive at a computable posterior distribution. We may use a measurement operator of the form
\begin{equation}
	R^{i j}_{\; \; t x} = M^i_{\; t} \ B^j_{\; x} = \delta(t_i - t) \ \delta(x_j - x) \ ,
\end{equation}
i.E. each measurement singles out a specific location in space-time. We notice that arbitrary (e.g. random) space-time locations again renders the simulation to scale with $N^2$ ($N^3$). To minimize this computational burden more sophisticated methods of choosing design points in space-time have been proposed. E.g.\ \cite{inproceedings} aims to choose design points such that the posterior uncertainty is minimized, i.E.\ by minimizing the trace or the determinant of the posterior covariance w.r.t.\ the locations of the design points.
For many PDEs, however, it is important to satisfy the equation at many points simultaneously in order to arrive at a good numerical approximation. Therefore, in this work, we make use of the specific prior structure to arrive at an almost linear scaling of the proposed method.

To this end we notice that due to the Markov property of the IWP, the distribution at a later time, given all Fourier modes in the past, only depends on the latest Fourier modes. In analogy to Eq.\ \eqref{eq:iwpprior}, for each Fourier mode $k$ we get an independent Markov process of the form
\begin{align}\label{eq:iwppriormodes}
&\left.P\left(\left(\begin{matrix}
\tilde{s}^{i k} \\ \dot{\tilde{s}}^{i k}
\end{matrix}\right) \right| \left(\begin{matrix}
\tilde{s}^{(i-1) k} \\ \dot{\tilde{s}}^{(i-1) k}
\end{matrix}\right) \right) \notag\\ &= \mathcal{G}\left(\left(\begin{matrix}
\tilde{s}^{i k} \\ \dot{\tilde{s}}^{i k}
\end{matrix}\right) - \left(\begin{matrix}
1 & \Delta_i \\
0 & 1
\end{matrix}\right) \left(\begin{matrix}
\tilde{s}^{(i-1) k} \\ \dot{\tilde{s}}^{(i-1) k}
\end{matrix}\right), \left|\sigma^k\right|^2 \left(\begin{matrix}
\nicefrac{\Delta_i^3}{3} & \nicefrac{\Delta_i^2}{2} \\
\nicefrac{\Delta_i^2}{2} & \Delta_i
\end{matrix}\right)\right),
\end{align}
with $\tilde{s}^{i k} = (M \tilde{s})^{i k}$.

However, the process only remains Markov if we keep all (infinitely many) modes in memory. If we additionally require the spatial locations to be on the same regular grid with $K$ points, i.E. $x_j = \nicefrac{j}{K}$, we notice that we can construct a discrete Markov process since
\begin{align}\label{eq:fouriersource}
	&e^{2 \pi i (k + n K) x_j} = e^{2 \pi i k} e^{2 \pi i n K \nicefrac{j}{K}} = e^{2 \pi i k} \\ &\forall j \in \left\lbrace 0, 1, ..., K-1 \right\rbrace \ , \ n \in \mathds{Z} \ . \notag
\end{align}
I.e.\ each Fourier mode $k$ shifted by multiples of $K$ coincides with the mode $k$ for each location on the grid.
Consequently we can represent the field values on the grid using only $K$ modes as
\begin{align}\label{eq:closure}
	\left(\bar{s}^{(c)}\right)^{t j} &\equiv \left(B s^{(c)}\right)^{t j} = \sum_{k = -\nicefrac{K}{2}+1}^{\nicefrac{K}{2}} \left(\tilde{\bar{s}}^{(c)}\right)^{t k} e^{2 \pi i k x_j} \notag\\ &\equiv \mathcal{F}^j_k \left(\tilde{\bar{s}}^{(c)}\right)^{t k} \ ,
\end{align}
where we defined the discrete Fourier transformation $\mathcal{F}$. The finite Fourier modes $\tilde{\bar{s}}$ are defined in terms of $\tilde{s}$ as
\begin{equation}
	\left(\tilde{\bar{s}}^{(c)}\right)^{t k} = \sum_{n = -\infty}^\infty \left(2 \pi i \left(k+n K\right)\right)^c \tilde{s}^{t (k + n K)} \ .
\end{equation}
Each discrete Fourier mode can be expressed in terms of an infinite sum of Gaussian random variables and thus itself is Gaussian. Note that for each spatial derivative $c$, however, the terms within the sum are different and therefore the summation results in a vector $\tilde{\bar{\mathbf{s}}} = \left(\tilde{\bar{s}}^{(0)}, \tilde{\bar{s}}^{(1)}, ...\right)$ of correlated Gaussian random variables, one for each spatial derivative involved in the PDE. The reason for this is that even though the field and its derivatives can be represented on the same grid, taking the derivative does not commute with the discretization operation $B$.

The infinite Fourier modes $\tilde{s}^k$ are solutions of the IWP process defined in Eq.\ \eqref{eq:iwppriormodes}, and therefore we may use an analogous derivation for the discrete representation of the time derivatives $\dot{ \tilde{\bar{\mathbf{s}}}} = \left(\dot{\tilde{\bar{s}}}^{(0)}, \dot{\tilde{\bar{s}}}^{(1)}, ...\right)$ to arrive at a discrete Markov prior of the form
\begin{align}\label{eq:pdemarkovprior}
	&\left.P\left(\left(\begin{matrix}
	\tilde{\bar{\mathbf{s}}}^{i k} \\ \dot{\tilde{\bar{\mathbf{s}}}}^{i k}
	\end{matrix}\right)\right|\left(\begin{matrix}
	\tilde{\bar{\mathbf{s}}}^{(i-1) k} \\ \dot{\tilde{\bar{\mathbf{s}}}}^{(i-1) k}
	\end{matrix}\right)\right) \notag\\ &= \mathcal{G}\left(\left(\begin{matrix}
	\tilde{\bar{\mathbf{s}}}^{i k} \\ \dot{\tilde{\bar{\mathbf{s}}}}^{i k}
	\end{matrix}\right) - \left(\begin{matrix}
	1 & \Delta_i \\ 0 & 1
	\end{matrix}\right)  \left(\begin{matrix}
	\tilde{\bar{\mathbf{s}}}^{(i-1) k} \\ \dot{\tilde{\bar{\mathbf{s}}}}^{(i-1) k}
	\end{matrix}\right) , \left(\begin{matrix}
	\nicefrac{\Delta_i^3}{3} & \nicefrac{\Delta_i^2}{2} \\ \nicefrac{\Delta_i^2}{2} & \Delta_i
	\end{matrix}\right) \otimes \mathbf{D}^k \right) \notag\\ & \forall k \in [-\nicefrac{K}{2} + 1, \nicefrac{K}{2}] \ , 
\end{align}
where $\tilde{\bar{\mathbf{s}}}^{i k} = M^i_t \tilde{\bar{\mathbf{s}}}^{t k}$ and $\otimes$ denotes the tensor product. The discrete Fourier mode covariance $\mathbf{D}^k$ takes the form
\begin{align}\label{eq:cov}
\left(\mathbf{D}^k\right)^{c d} &\equiv \left<\left(\tilde{\bar{s}}^{(c)}\right)^k \left(\tilde{\bar{s}}^{(d)}\right)^k\right> \notag\\ &= (-1)^d \sum_{n = -\infty}^\infty  \left(2 \pi i (k+nK)\right)^{c+d} \left|\sigma^{k+nK}\right|^2 \ .
\end{align}
The Markov property of the IWP remains in the discrete representation of the field since we defined the space and time correlations to be independent a priori.
See Appendix \ref{ap:distribution} for a derivation of $\mathbf{D}^k$.

The discrete Fourier transformation defined in Eq.\ \eqref{eq:closure} is invertible, and therefore we can construct the measurement equation associated with the PDE (Eq.\ \eqref{eq:pde}) in terms of the Fourier modes as
\begin{equation}
	d^{i k} = \dot{\tilde{\bar{s}}}^{i k} - \left(\mathcal{F}^{-1} f\left(\mathcal{F} \tilde{\bar{s}}^{(0)}, \mathcal{F} \tilde{\bar{s}}^{(1)}, ...\right)\right)^{i k} \equiv \left(\dot{\tilde{\bar{s}}} - g\left(\tilde{\bar{\mathbf{s}}}\right)\right)^{i k} \ .
\end{equation}

\subsubsection{Posterior distribution}
In direct analogy to the ODE setting, we can combine the observational data $d$ with the prior to construct a posterior distribution. Let $\mathbf{u} = \left(\tilde{\bar{s}}^{(0)}, \tilde{\bar{s}}^{(1)}, ...\right)$ be the discretized Fourier space field values and their higher order spatial derivatives and $\mathbf{v} \equiv \left(\dot{\tilde{\bar{s}}}^{(1)}, ...\right)$ be the time derivative of the spatial derivatives in $\mathbf{u}$, we get that
\begin{align}\label{eq:pdeposterior}
	&P\left(\mathbf{u}, \mathbf{v} | d = 0, \mathbf{u}^0, \mathbf{v}^0\right) \notag\\ &\propto \prod_{i=1}^{N-1} P\left(\left(\begin{matrix}
	\mathbf{u}^i \\
	\dot{\tilde{\bar{s}}}^i = g\left(\mathbf{u}^i\right) \\
	\mathbf{v}^i
	\end{matrix}\right) \left| \left(\begin{matrix}
	\mathbf{u}^{i-1} \\
	\dot{\tilde{\bar{s}}}^{i-1} = g\left(\mathbf{u}^{i-1}\right) \\
	\mathbf{v}^{i-1}
	\end{matrix}\right)\right)\right. \notag\\ &= \prod_{i=1}^{N-1}
	\left[P\left(\mathbf{v}^i \left| \left(\begin{matrix}
	\mathbf{u}^i \\
	\dot{\tilde{\bar{s}}}^i = g\left(\mathbf{u}^i\right)
	\end{matrix}\right) ,
	\left(\begin{matrix}
	\mathbf{u}^{i-1} \\
	\dot{\tilde{\bar{s}}}^{i-1} = g\left(\mathbf{u}^{i-1}\right) \\
	\mathbf{v}^{i-1}
	\end{matrix}\right)
	\right)\right.\right. \notag\\ &\left.P\left( \left(\begin{matrix}
	\mathbf{u}^i \\
	\dot{\tilde{\bar{s}}}^i = g\left(\mathbf{u}^i\right)
	\end{matrix}\right) \left| \left(\begin{matrix}
	\mathbf{u}^{i-1} \\
	\dot{\tilde{\bar{s}}}^{i-1} = g\left(\mathbf{u}^{i-1}\right)\\
	\mathbf{v}^{i-1}
	\end{matrix}\right) \right)\right.\right] \ .
\end{align}
Here, the involved conditional distributions can be directly constructed from Eq.\ \eqref{eq:pdemarkovprior}. We notice that the distribution of $\mathbf{v}^i$ remains Gaussian and we can directly sample it once we solved the simulation step for $\mathbf{u}^i$ by constructing the conditional distribution of $\mathbf{v}^i$ from Eq.\ \eqref{eq:pdemarkovprior}. The distribution of $\mathbf{u}^i$ may again be rewritten in terms of a non-linear filter as
\begin{align}\label{eq:pdepartpos}
	&P\left( \left(\begin{matrix}
	\mathbf{u}^i \\
	\dot{\tilde{\bar{s}}}^i = g\left(\mathbf{u}^i\right)
	\end{matrix}\right) \left| \left(\begin{matrix}
	\mathbf{u}^{i-1} \\
	\dot{\tilde{\bar{s}}}^{i-1} = g\left(\mathbf{u}^{i-1}\right)\\
	\mathbf{v}^{i-1}
	\end{matrix}\right) \right)\right. = \notag\\ &= P\left(\dot{\tilde{\bar{s}}}^{i} = g\left(\mathbf{u}^i\right) \left| \mathbf{u}^i , \left(\begin{matrix}
	\mathbf{u}^{i-1} \\
	\dot{\tilde{\bar{s}}}^{i-1} = g\left(\mathbf{u}^{i-1}\right)\\
	\mathbf{v}^{i-1}
	\end{matrix}\right)\right)\right. \times \notag\\ &P\left(\mathbf{u}^i \left| \left(\begin{matrix}
	\mathbf{u}^{i-1} \\
	\dot{\tilde{\bar{s}}}^{i-1} = g\left(\mathbf{u}^{i-1}\right)\\
	\mathbf{v}^{i-1}
	\end{matrix}\right)\right)\right.
\end{align}

Eq.\ \eqref{eq:pdeposterior} and \eqref{eq:pdepartpos} describe the central results of our work. Under the given prior assumptions and measurement setting the posterior becomes a Markov process in time in the finite state vector $\mathbf{u}^i$. Furthermore, each time step is presented as a non-linear Bayesian filtering problem, where the second probability on the r.h.s.\ in Eq.\ \eqref{eq:pdepartpos} is a Gaussian prior distribution in $\mathbf{u}^i$ that acts as a predictive step to construct the next step from the previous one. The first distribution may be regarded as a (in general non-linear) likelihood which acts as a regularization by comparing the time derivative $\dot{\tilde{\bar{s}}}^i$ constructed via the PDE from $\mathbf{u}^i$, to the conditional distribution of $\dot{\tilde{\bar{s}}}^i$ that arises from the previous step and the prior process. See Algorithm \ref{al:fixedpspec} for a pseudo-code description of the resulting algorithm.

\subsection{Posterior properties}
It is noteworthy that, in contrast to the ODE setting, even though we use an IWP prior in time, it is in general not sufficient to only store the field values on the grid. We also have to keep the involved spatial derivatives $\mathbf{u}$ and, maybe even more surprising, the spatial derivatives of the first time derivative $\mathbf{v}$ in memory, in order to be fully consistent with the continuous prior process. In fact, as the spatial derivatives of the first time derivative  do not enter the PDE, we may analytically integrate over these quantities, but the resulting process would loose the Markov property, which we believe is in general not desirable. However, as we have seen, once we have solved the inference problem for $\mathbf{u}$ we can directly sample $\mathbf{v}$ as the conditional distribution remains Gaussian.

On the other hand, given a fixed step size, the spatial resolution, and a spectrum $|\sigma|^2$, we may rewrite the posterior distribution in terms of the generative process associated with the predictive prior of $\mathbf{u}$. This reads
\begin{equation}
	\mathbf{u}^{i} = \mathbf{u}^{i-1} + \Delta_i \left(\begin{matrix}
	g\left(\mathbf{u}^{i-1}\right) \\
	\mathbf{v}^{i-1}
	\end{matrix}\right) + \sqrt{\nicefrac{\Delta_i^3}{3}} \ \mathbf{U} \Lambda \mathbf{r}^i \ ,
\end{equation}
with $\mathbf{r}^i \sim \mathcal{G}(\mathbf{r}^i, \mathds{1})$, and where $\mathbf{U} \Lambda \mathbf{U}^\dagger$ denotes the eigen-decomposition of the prior covariance $\mathbf{D}$ with $\mathbf{U}$ being a unitary matrix and $\Lambda$ a real diagonal matrix. Note that due to the homogeneity of the prior this covariance takes a block diagonal form in $k$ and therefore we only need to decompose a set of $K$ independent $(o+1)$-dimensional matrices where $o$ is the highest spatial derivative involved in the PDE. We notice that for fast decaying spectra (resulting in a strong spatial smoothness) the eigenvalues $\Lambda$ also decrease very fast. That means that we can define a precision level prior to the simulation up to which we want to keep track of discretization contributions, and set all eigenvalues below this threshold and all associated components in $\mathbf{r}^i$ to zero. This may reduce the burden of storing additional quantities on the grid.

\subsection{Power spectrum estimation}
So far we only considered the case of a given prior power spectrum $\sigma$. In practical applications, however, it is often unclear prior to the simulation which spatial correlation structure one should choose given the initial state and the PDE. A strongly decaying spectrum enforcing too much smoothness might result in a poor performance of the simulation algorithm as small scale structures are missing while a very flat spectrum might over-represent these scales and consequently leads to very high uncertainties.

Since we have formulated the simulation problem by means of Bayesian inference, it is straightforward to elevate the power spectrum to an unknown quantity that has to be inferred along with the solution. To this end we may write
\begin{equation}\label{eq:pspecpost}
	P(\mathbf{u}, \mathbf{v}, \sigma | d, \mathbf{u}^0, \mathbf{v}^0) \propto P(\mathbf{u}, \mathbf{v} | d , \mathbf{u}^0, \mathbf{v}^0, \sigma) \ P(\sigma) \ ,
\end{equation}
where $P(\mathbf{u}, \mathbf{v} | d , \mathbf{u}^0, \mathbf{v}^0, \sigma)$ is defined via Eq.\ \eqref{eq:pdeposterior}.

A more difficult question is how to construct a useful prior distribution for $\sigma$ as in order to construct the distribution of $\mathbf{u}$ and $\mathbf{v}$ we have to compute the infinite sums associated with $\mathbf{D}$ (see Eq.\ \eqref{eq:cov}). In this work we follow an approach originally developed for power spectrum estimation within the context of Bayesian imaging \cite{refId0}: First consider the spectrum on a double logarithmic scale as
\begin{equation}
	\sigma(|k|) = e^{\tau(l)} \quad \text{with} \quad l = \log(|k|) \ .
\end{equation}
This provides a useful scale for power spectra as power laws appear as straight lines on this scale. As power-law shaped spectra are reasonable for many physical processes, we aim to construct a prior that, in absence of further information, follows a power law. Furthermore we require that deviations from this power-law are smooth (i.E. differentiable) on log-log-scale.
To this end we assume that $\tau$ solves an IWP process in the log-coordinates $l$ of the form
\begin{equation}\label{eq:tauprocess}
	\frac{\partial^2 \tau}{\partial l^2} = \sigma_\tau \xi \quad \text{with} \quad  \xi \sim \mathcal{G}(\xi, \mathds{1}) \ ,
\end{equation}
where $\sigma_\tau$ is a positive scaling factor.
Finally, we realize this process on a regular grid in $l$ with $L$ pixels, up to a maximal value $l_{\mathrm{max}}$, and approximate all intermediate values of $\tau$ via bi-linear interpolation in $l$. This allows us to approximately compute the covariance $\mathbf{D}$ by summing up all contributions to the sum up to $n_{\mathrm{max}}$ with $l_{\mathrm{max}} = \log\left(n_{\mathrm{max}} K \right)$ where $K$ is the number of pixels of the spatial grid. The bi-linear interpolation additionally allows to approximately compute the sum directly from the values of $\tau$ on the logarithmic grid $l$ without the need to realize a high resolved version of $\sigma$ on linear scale. Furthermore, as we define a regular grid on logarithmic scale in $|k|$ we can easily extend the spectrum to extremely large values of $|k|$ (large $n_{\mathrm{max}}$), far below the smallest resolved scales of the simulation. For a detailed discussion of these prior properties see e.g. \cite{refId0} and \cite{arras2020variable}.

We notice that a time invariant spectrum constructed this way renders the full posterior to be non-Markov since all steps depend on the same spectrum. We can restore the Markov property by introducing a different spectrum for each time step $\tau^i$. Specifically, we assume the spectrum to be piecewise constant for the length of the time step, but different for each step.
Furthermore, to increase stability, we may assume that the power spectra of subsequent steps are correlated, which is a reasonable assumption since we do not expect the statistical properties to vary arbitrarily strong between two subsequent time steps. A simple way to introduce such correlations is by assuming that $\tau$ follows a discrete time Wiener process, that is
\begin{equation}
	\tau^i = \tau^{i-1} + \Delta_i \tilde{\tau}^i \ .
\end{equation}
Specifically the current log-spectrum $\tau^i$ can be constructed from the previous one $\tau^{i-1}$ and a random component $\tilde{\tau}^i$. We let $\tilde{\tau}^i$ be distributed according to an IWP in the log-Fourier coordinates $l$, as defined via Eq.\ \eqref{eq:tauprocess}. This renders the full time-Fourier process for $\tau$ to be a discrete Wiener Process in time and an IWP in the log-Fourier coordinates $l$.

\subsection{Composed algorithm}
The full algorithm using power spectrum estimation may be denoted as:

Given the previous state $X^{(i-1)} = \left(\mathbf{u}{(i-1)}, \mathbf{v}{(i-1)}, \tau^{(i-1)}\right)$, use the posterior distribution constructed from Eq.\ \eqref{eq:pdepartpos} and Eq.\ \eqref{eq:pspecpost} to compute an estimate (or sample) for $\mathbf{u}^i$ and $\tau^i$ via e.g.\ a joint Maximum a Posteriori (MAP) estimate, a Variational approximation, or Monte Carlo based sampling. Use this estimate (sample) in the distribution of $\mathbf{v}^i$ (see Eq.\ \eqref{eq:pdeposterior}) to sample $\mathbf{v}^i$ conditional to $\mathbf{u}^i$, $\tau^i$ and the previous state $X^{(i-1)}$. Given the new full state $X^i$ we may repeat the procedure to compose a new time-step. For a pseudo code representation see Algorithm \ref{al:variablepspec}

\begin{figure}[htp]
	\begin{algorithm}[H]
		\caption{PDE simulation with fixed spectrum}
		\label{al:fixedpspec}
		\begin{algorithmic}[]
			\State {Input: $\mathbf{u}^0$, $\mathbf{v}^0$, $\sigma$, PDE}
			\For{$i = 1$ to $N$}
			\State {Given $\left(\mathbf{u}^{i-1}, \mathbf{v}^{i-1}\right)$ and $\sigma$, solve Bayesian filtering problem (Eq.\ \eqref{eq:pdepartpos}) to get an estimate (sample) for $\mathbf{u}^i$}
			
			\State {Given $\mathbf{u}^i$ use Eq.\ \eqref{eq:pdeposterior} to sample $\mathbf{v}^i$}
			\EndFor
			
			\Return $\left\lbrace \left(\mathbf{u}^i, \mathbf{v}^i\right)\right\rbrace_{i \in \lbrace 1, ..., N\rbrace}$
		\end{algorithmic}
		
	\end{algorithm}
\end{figure}

\begin{figure}[htp]
	\begin{algorithm}[H]
		\caption{PDE simulation with variable spectrum}
		\label{al:variablepspec}
		\begin{algorithmic}[]
			\State {Input: $\mathbf{X}^0$, PDE}
			\For{$i = 1$ to $N$}
			\State {Given $\mathbf{X}^{i-1}$, solve the joint Bayesian filtering problem of Eqs.\ \eqref{eq:pdepartpos} and \eqref{eq:pspecpost} to get an estimate (sample) for $\mathbf{u}^i$ and $\tau^i$}
			
			\State {Given $\mathbf{u}^i$ and $\tau^i$ use Eq.\ \eqref{eq:pdeposterior} to sample $\mathbf{v}^i$}
			
			\State Set $\mathbf{X}^i = \left(\mathbf{u}^i, \mathbf{v}^i, \tau^i\right)$
			\EndFor
			
			\Return $\left\lbrace \mathbf{X}^i \right\rbrace_{i \in \lbrace 1, ..., N\rbrace}$
		\end{algorithmic}
		
	\end{algorithm}
\end{figure}

\subsubsection{Initial conditions}
We notice that initial conditions $s^0$, evaluated on the grid, do not fully determine the initial state $X^0$ that is needed to start the simulation as $X^0$ also consists of the spatial derivatives of the continous field, evaluated on the grid, and the initial power spectrum $\tau^0$.
However, there are multiple ways to estimate an initial state $X^0$ given $s^0$. For example we may estimate the large scale (scales that are resolved by the simulation grid) power spectrum from the initial conditions directly and accompany this estimate with a consistent initial guess for the small scale spectrum. Given this spectrum, it is straightforward to estimate the spatial derivatives needed for $X^0$, given the spectrum and $s^0$ via Gaussian regression. We may even perform a probabilistic estimate and sample from the corresponding distribution to construct $X^0$ in order to propagate the uncertainty that arises from insufficient knowledge of the initial state into the simulation.

In this work, however, we want to study the performance of the simulation algorithm itself, and therefore assume that the initial state $X^0$ is fully given, i.E.\ we start with an initial condition that allows us to compute the spatial derivatives analytically.

\section{Applications}
In the following we present the application of the proposed methods to two systems, the diffusion equation as well as the viscous Burgers equation. All applications are conducted on the same regular grid in space, with $128$ pixels and periodic boundary conditions. The power spectra are realized on a logarithmic regular grid with $500$ pixels and a maximal value $l_{\mathrm{max}}$ corresponding to an effective Fourier space $100$ times the resolution of the simulation grid. This large effective Fourier space ensures that, at any point in the given examples, the spectra are numerically zero outside this region.

\subsection{Diffusion equation}
\begin{figure*}[htp]
	\centering
	\includegraphics[scale=.65, angle=0]{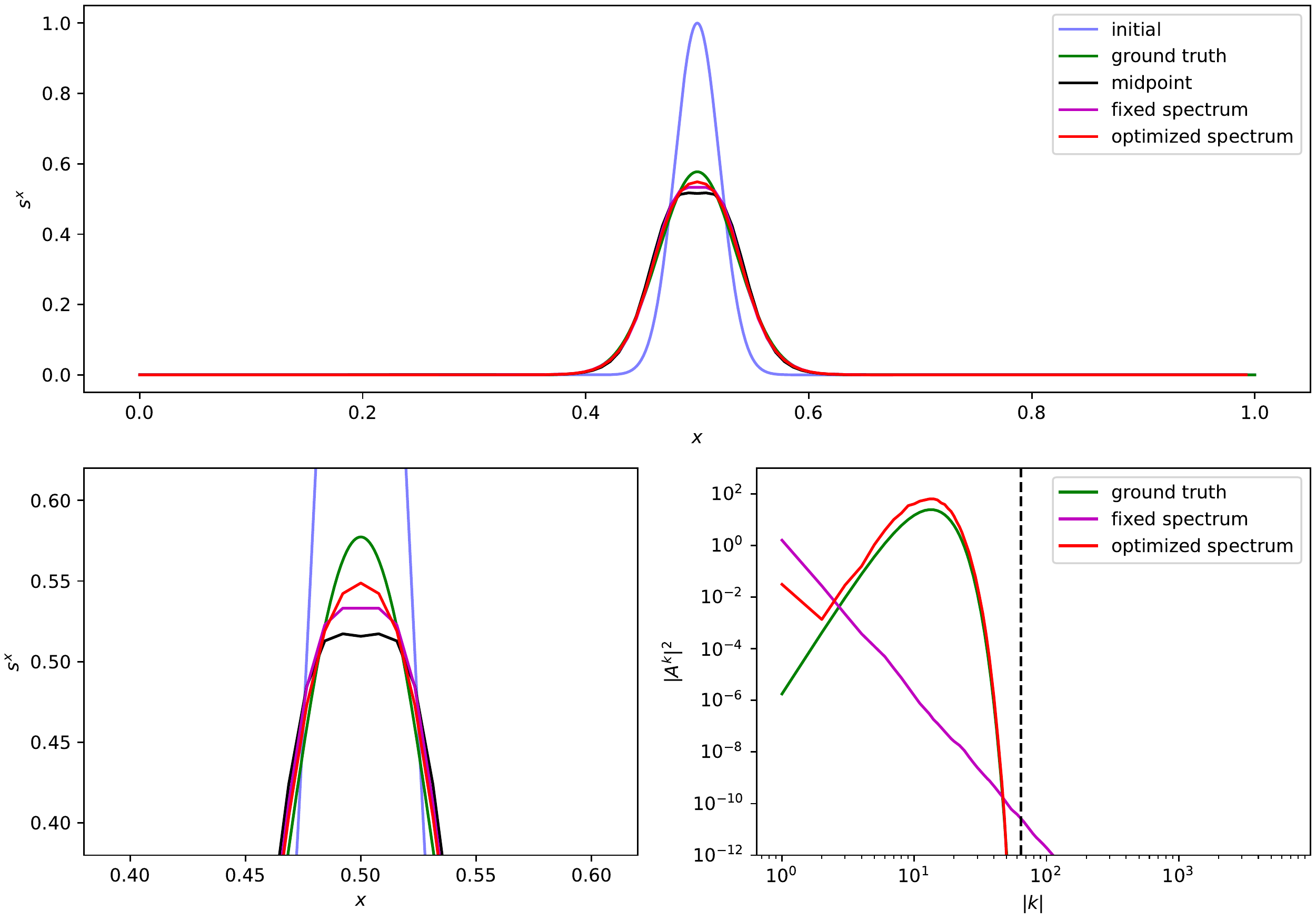}
	\centering
	\caption{{\bf Top:} First time step of the simulation of the diffusion equation with an initial Gaussian profile (blue). The green line corresponds to the ground truth, the black line to the midpoint rule, the purple line to the posterior mean of the reconstruction using a fixed power spectrum $\propto |k|^{-6}$, and the red line corresponds to the MAP estimate of the simulation with an adaptive power spectrum. {\bf Bottom left:} Detailed version of the simulation step zoomed into the central region. {\bf Bottom right:} Power spectra of the simulation on a double-logarithmic scale. Purple: Spectrum of the simulation step with a fixed spectrum. Red: MAP estimate of the optimized spectrum. Green: Ground truth of the spectrum. Here ground truth refers to the spectrum that was reconstructed using the true time evolution as a realization of the corresponding Gaussian prior distribution. The black dashed line indicates the largest harmonic mode corresponding to the resolution of the simulation.} \label{fig:diffusion_step} 

	\centering
	\includegraphics[scale=.65, angle=0]{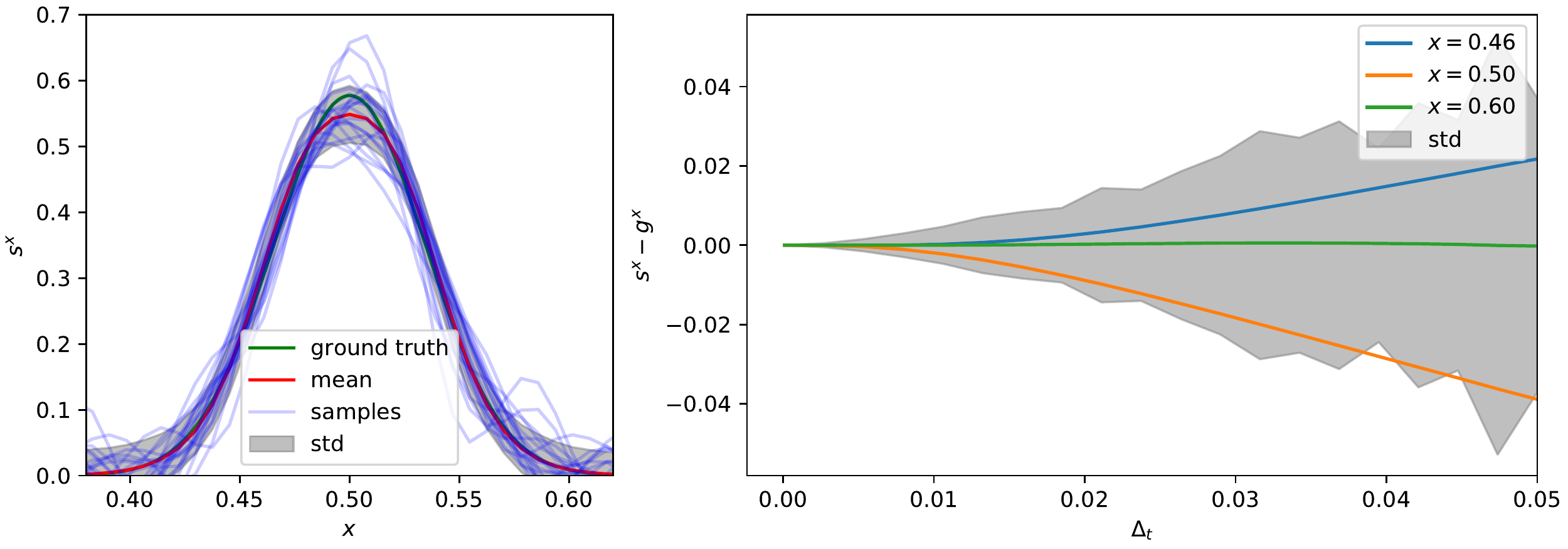}
	\centering
	\caption{{\bf Left:} Ground truth (green), posterior mean (red), posterior samples (light blue), and posterior standard deviation (gray) of the first time step of the diffusion equation. The posterior samples as well as the standard deviation were conducted by means of the empirical Bayes' approach. Specifically, the posterior distribution conditional to the MAP estimate of the optimized spectrum is used.
	{\bf Right:} Colored lines: Residual difference between the ground truth and the posterior mean at multiple locations of the spatial domain as a function of step size $\Delta_t$. The corresponding posterior standard deviation (valid for any location) is given as the gray contour.} \label{fig:diffusion_unc} 
\end{figure*}

To emphasize the influence of the spectrum on the simulation we start with the simple case of a diffusion equation, that is
\begin{equation}
	\dot{s} = f(s) = \nu \ s^{(2)} \ , \ \nu > 0 \ ,
\end{equation}
and choose a Gaussian profile as the initial state. In Figure \ref{fig:diffusion_step} we depict the MAP estimate of the first step for a step size of $\Delta_1 = 0.04$, and for $\nu = 0.01$. We show two different modes of the simulation scheme: the case of a given generic power spectrum of the form $|\sigma^k|^2 \propto |k|^{-6}$ as well as the case where we optimize for the spectrum together with the solution. As a comparison, we also compute the solution given by the trapezodial rule, where in this case the spatial derivatives are computed via discrete Fourier derivatives, i.E.\ $\left(s^{(2)}\right)^k = (2 \pi i k)^2 \left(s^{(0)}\right)^k$. This method may serve as a standard comparison as it also requires the differential equation to be satisfied for the current as well as the future state simultaneously and therefore is an implicit method of second order, such as the two approaches proposed in this work are. We see in Figure \ref{fig:diffusion_step} that compared to the standard method, both approaches are closer to the ground truth, with the optimized spectrum being slightly closer.

Furthermore, in Figure \ref{fig:diffusion_unc}, we compare the ground truth to the posterior mean of the simulation and also depict the posterior uncertainty of the problem. We approximate the posterior distribution via the empirical Bayes approach, that is, we use the Maximum a posterior (MAP) estimate of the logarithmic power spectrum $\tau^*$ and compute the conditional posterior distribution of the solution $s$, given $\tau^*$. This conditional posterior is analytically computable since the linear dynamics together with a Gaussian prior distribution results in a Gaussian posterior for $s$, given $\tau^*$. We see that the posterior mean is in agreement with the ground truth within posterior uncertainties. Furthermore, on the right hand side of Figure \ref{fig:diffusion_unc}, we depict the residual between the ground truth and the reconstruction as a function of the step size for various locations. Again, the deviation agrees with the uncertainties and furthermore we notice that due to the fact that the prior is stationary, and the diffusion equation is linear and stationary, the posterior distribution also remains a stationary process in space and therefore the posterior uncertainty is the same for every location.

\begin{figure*}[htp]
	\centering
	\includegraphics[scale=.65, angle=0]{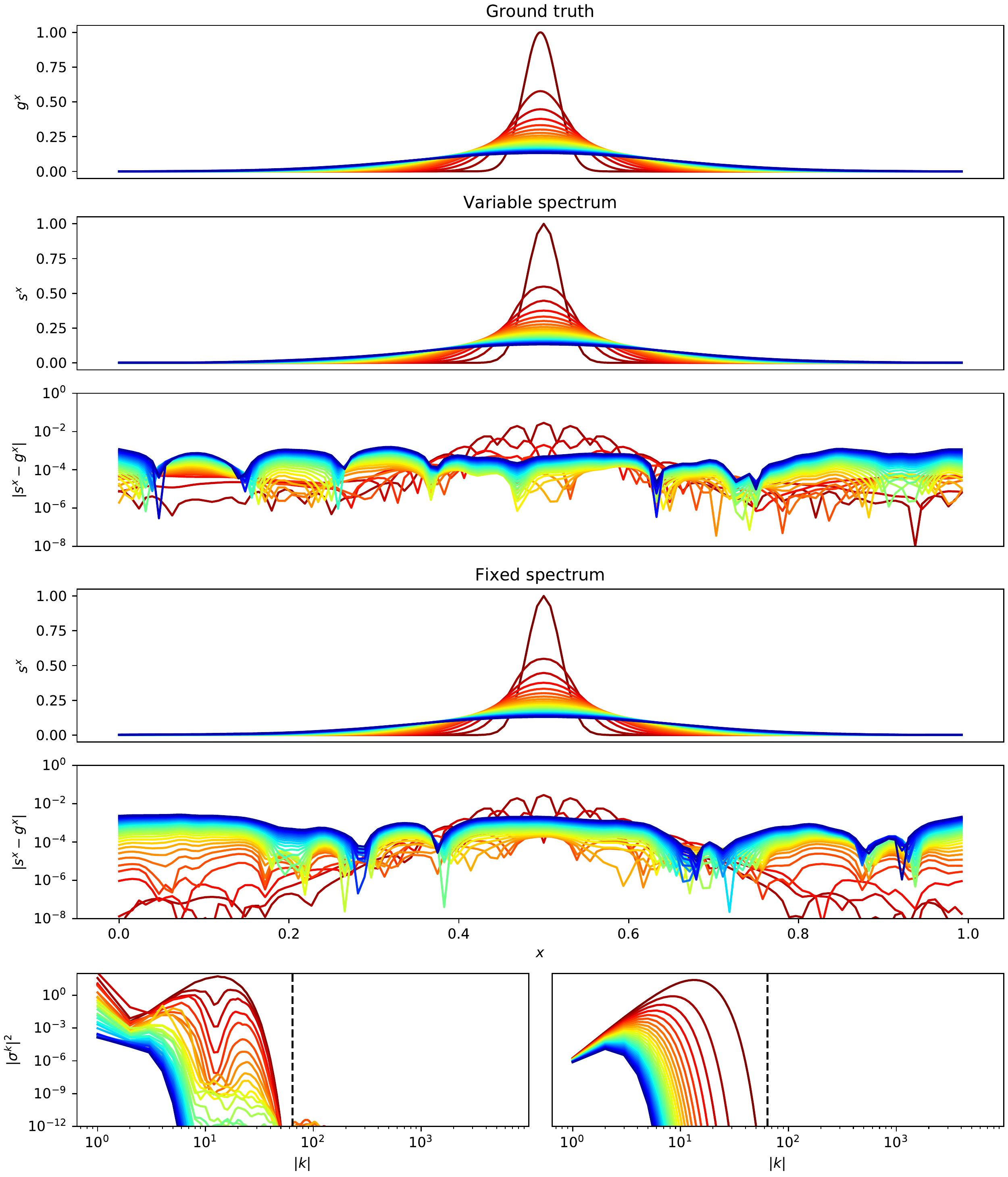}
	\centering
	\caption{Color coded time evolution of the diffusion equation. Red indicates early times and blue indicates the latest time-steps. {\bf Top to bottom:} Ground truth, reconstruction using a variable spectrum (i.E.\ joint optimization for solution and spectrum), residual norm between ground truth and reconstruction, reconstruction using the fixed spectrum derived from the ground truth, and corresponding residual norm.
	{\bf Bottom left:} Reconstructed power spectra for each time-step of the joint optimization case.
	{\bf Bottom Right:} Power spectra computed from the ground truth.} \label{fig:diffusion_time} 
\end{figure*}

Finally, in Figure \ref{fig:diffusion_time} we depict the time evolution of the simulation together with the ground truth and the estimated power spectra for every time step. As a comparison, we also depict the time evolution for a simulation setting where we used the power spectra computed from subsequent steps of the ground truth, and solved the simulation problem conditional to these spectra.

We see that as time progresses, the initially sharp spatial distribution tends to decay and smooth out over the spatial domain. Consequently, the reconstructed power spectra show less power on small scales as time progresses and only large scale power remains. Furthermore, the overall magnitude of the power spectrum decreases, which indicates that the uncertainty (and therefore the local error) of later time steps become smaller. This adaptive control of the spectrum leads to a better quantification of the local error and therefore also leads to a more sophisticated control of the global error of the system. We notice, however, that the inferred power spectra of intermediate steps are substantially different from the power spectra of the ground truth. First, on the largest scales the reconstructed power spectra has more power compared to the ground truth. This is a common issue that appears when jointly inferring a field with its power spectrum, as for these modes inference is very degenerate and consequently mostly dominated by the prior assumptions. A more suitable prior in terms of more restrictive hyper-parameters might improve this behaveiour. The second difference becomes apparent for small scale modes where there is too much power around $|k| \in  [20, 30]$. We believe that this effect is rooted in the large step size of the given simulation setting: The first steps of the ground truth show a rapid decay of these modes which cannot fully be captured by the simulation step and thus power remains on these scales that gets picked up by the power spectra estimate. However, as time progresses, the power of these scales eventually decay due to the diffusive dynamics of the process.

\subsection{Burger's equation}
As a second example, we study the performance of the proposed approach in the context of the (viscous) Burgers equation. Specifically
\begin{equation}
	\dot{s} + \ s \ s^{(1)} = \nu \ s^{(2)} \ .
\end{equation}
We again start with a Gaussian profile as the initial state and set $\Delta_i = 3 \times 10^{-3}$ and $\nu = 4 \times 10^{-3}$. 

\begin{figure*}[htp]
	\centering
	\includegraphics[scale=.65, angle=0]{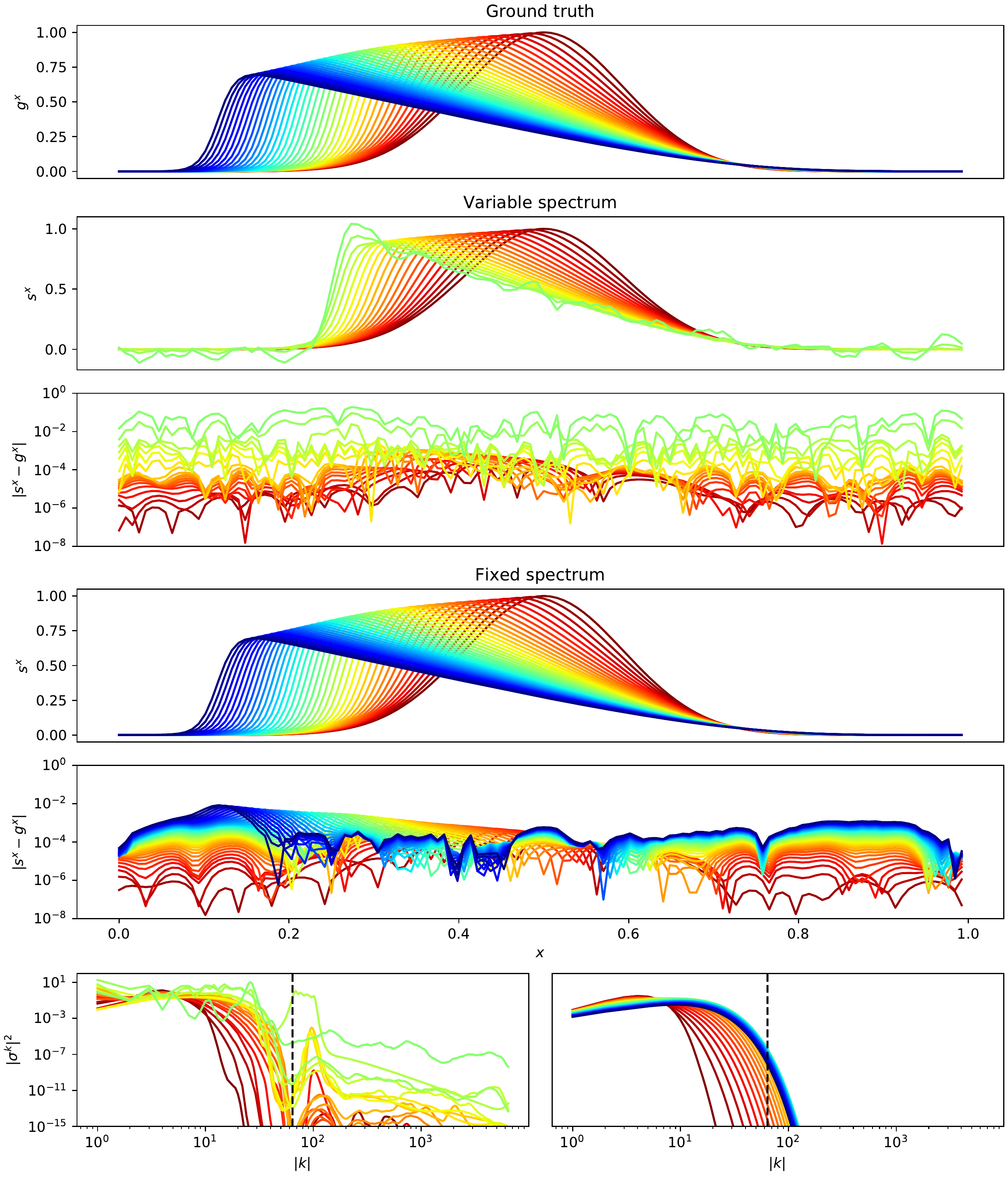}
	\centering
	\caption{Same composition as Figure \ref{fig:diffusion_time}, but for the time evolution of the Burgers equation.} \label{fig:burger_time} 
\end{figure*}

The Burger's equation is known to develop strong shock waves for small viscosity $\nu$, which means that in contrast to the diffusion equation, small scale structures become more relevant as time progresses. Indeed we find that if we compute the power spectra of subsequent time steps from the ground truth (see bottom right of Figure \ref{fig:burger_time}) we see how the spectrum gains power on small scales, while the large scale power remains almost unchanged. In addition we also notice that after a few time steps there is non-negligible power on scales that are smaller then the smallest resolved scales of the simulation grid.

It turns out that, when applying the adaptive simulation to this setup (see Figure \ref{fig:burger_time}), it is only possible to consistently infer the power spectra along with the solution for scales that are also resolved by the simulation grid. As we only require the differential equation to be satisfied on the grid, there is no direct information about smaller scales that enter the reconstruction and therefore the power spectrum estimation, and ultimately also the simulation itself breaks down as the shock forms. This leads us to the conclusion that using only the feedback of small scales to the large scales provides insufficient information to properly infer the small scale statistics. Without further prior information, we believe that the only way to properly access these scales is via resolving them on a grid with high enough resolution.

However, we notice that it is possible to circumvent the need of realizing the process on a high resolution grid, via the usage of appropriate prior information. To this end consider the middle panels of Figure \ref{fig:burger_time}, where we used the power spectra estimated from the ground truth to construct a simulation scheme with fixed spectrum on the same resolution as the adaptive one (i.E. a spatial discretization of $128$ pixels). It turns out that in contrast to the adaptive scheme, the simulation remains stable and is in agreement with the ground truth long after the adaptive scheme diverged. This result highlights the second key mechanism of a probabilistic treatment of PDE simulation: even though the spatial resolution appears to be insufficient to fully resolve the state, the consistent treatment of discretization via the introduction of spatial derivatives as additional random variables allows for a simulation that remains in agreement with the ground truth. As the correct power spectra are given in this setup, they provide small scale structures consistent with the given PDE and in turn allow for a correct feedback of the small (unresolved) scales to larger (resolved) scales.

\section{Comparison to IFD}
In this work, as well as in IFD there exists the concept of a measurement operator $R$ that specifies the evaluated values of the field. In IFD the resulting measurements are the quantities that are ultimately stored on a computer for a given time-step, meaning that if $R$ singles out a finite set of spatial locations, as used in this work, the corresponding field values are stored. In contrast, in this work not only the field values but also the spatial derivatives involved in the PDE are stored. However, we note that one can alter the measurement operator of IFD to measure not only the values but also the spatial derivatives, to result at the same set of quantities that have to be stored. The important difference is that while in IFD this is a choice made by the user, in this work it is a result of the method in order to arrive at a computable distribution that is consistent with the continuous prior process.

Furthermore, in this work $R$ also defines the set of space-time locations at which the process has to fulfill the PDE. This is fundamentally different from IFD as IFD aims to fulfill the PDE at every location. As a consequence there is no need for a prior time correlation in IFD as, in case of a Gaussian prior, the only quantity necessary to translate between the finite state and the distribution of the field is a prior spatial correlation structure. However, we note that for most non-linear applications, the exact time evolution that is required for IFD is not available and thus an approximation has to be made there, which is not captured in a probabilistic fashion. Consequently uncertainties arising from approximated time evolution are not captured within IFD, while the approach in this work takes into account these uncertainties and aims to fill the time gap via the assumed prior time correlation structure. However, requiring the PDE to be satisfied only at a discrete set of locations is also problematic as we have seen, in particular when we aim to infer the prior correlation structure (i.E.\ the power spectra) on scales that are not resolved by $R$.

\section{Conclusion}
In this work we derived a fully probabilistic framework for simulation of PDEs subject to periodic boundary conditions. The proposed method makes use of continuous space-time Markov process priors that are stationary in space, and incorporates artificial observational constraints that require the PDE to be satisfied on a regular grid. The Markov property allows for a formulation of the posterior such that the distribution of the current state is only conditional on the state at the previous time-step. The state of the system, however, not only consists of the field values realized on the grid, but also consists of the values of all spatial derivatives involved in the PDE. Only if these random variables are kept track of, the discrete Markov realization is consistent with the continuous process. Furthermore, the usage of prior distributions that are stationary in space, together with sampling on a regular and periodic grid with $K$ pixels, allows for an efficient $K \log(K)$ scaling of a single step of the algorithm via incorporation of Fast Fourier Transforms.

The Bayesian analysis of the problem allows for inference of hyper parameters, such as the spatial correlation structure, i.E. the prior power spectrum, alongside with the solution of the simulation. To this end we incorporate a non-parametric method of power spectra estimation, originally developed for Bayesian imaging by means of information field theory. The resulting joint estimation of spectrum and realization of the process leads to a simulation scheme that is closer to the ground truth compared to a method with a fixed, generic spectrum, and also allows for a more sophisticated error analysis in terms of the posterior uncertainty. We notice, however, that without further prior information about the small scale statistics, the inference of the power spectrum is only valid up to scales that are resolved by the simulation grid. As we have seen in the application to the Burgers equation, once scales below the grid resolution become relevant for the solution, the estimation of the spectrum becomes inaccurate, and as a consequence the simulation starts to diverge from the true solution.
If an accurate estimation of the small scale spectra are available, however, we notice that it is possible to use these spectra for a low-resolution simulation that remains consistent with the high-resolution setting.

Finally we may conclude that the approach for probabilistic PDE simulation provides novel insights into the interplay between prior assumptions entering a simulation algorithm and the involved PDE. However, additional work, in particular concerning small (unresolved) scale statistics, has to be done in order to improve the performance and stability of the proposed approach.

On the other hand, in addition to Bayesian uncertainty quantification, a fully probabilistic approach to simulation enables several novel key properties compared to traditional numerical simulation. For example, as the analysis gives rise to a posterior probability distribution that may be separated into a generative prior and a likelihood, it is straightforward to incorporate the simulation into a larger inference framework, in order to estimate for example parameters of the PDE or initial conditions, from observational data.

In addition, modern day machine learning techniques can be used to speed up the simulation algorithm. In particular neural networks have already successfully been applied to simulation using training data composed via traditional numerical simulation as an input (see e.g.\ \cite{SIRIGNANO20181339}). On the other hand, to circumvent the need of generating training data, which might be very expensive, \cite{raissi2017physics} has demonstrated that is possible to train a neural network to approximate the solution directly by minimizing the squared norm of the deviations of the PDE from zero at a discrete set of space-time locations using only the initial state and the PDE as an input. However, in \cite{raissi2017physics}, it has also been demonstrated that training a network to reproduce the internal stages of a high-order Runge-Kutta scheme rather than solely minimizing the squared norm associated with the PDE, appears to be more efficient due to the additional prior assumptions incorporated in the Runge-Kutta scheme. As Runge-Kutta type methods have a probabilistic interpretation in terms of a Gaussian process prior \cite{SchoberDH2014}, these results indicate that on one hand, neural networks are capable of approximating simulation steps, and on the other hand that a probabilistic posterior distribution for simulation, as derived in this work, may provide a more sophisticated measure for neural-network training. Specifically the posterior distribution is informed about both, the differential equation being satisfied, and a notion of continuity (and differentiability) in space and time in terms of the prior assumptions.

All in all, we believe that the probabilistic approach to simulation, in particular in terms of probabilistic numerics, is capable to provide further insights into numerical simulation, and to generalize existing algorithms. However, further work has to be done in order to arrive at a class of simulation algorithms that are capable of tackling broader classes of physically relevant PDEs.

\bibliographystyle{apsrev4-1}
\bibliography{simulation.bib}

\section*{Acknowledgments}
We would like to thank Reimar Leike and Philipp Arras for fruitful discussions and constructive feedback throughout the development process.

\newpage

\appendix
\section{Discrete prior}\label{ap:distribution}
Consider a Gaussian random field $s^{t x}$ with $x \in [0,1]$ on a periodic domain and $t \in [t_0, \infty)$. Furthermore $s$ has statistically homogeneous and isotropic statistics in space and follows an IWP in time. Specifically:
\begin{align}
	s^{t x} &= \sum_{k = -\infty}^\infty \tilde{s}^{t k} e^{2 \pi i k x} \\ \label{eq:apiwp}
	\ddot{\tilde{s}}^{t k} &= \sigma^k \ \xi^{t k} \quad \text{with} \quad \xi \sim \mathcal{G}(\xi, \mathds{1}) \ . \\\label{eq:apiwp2}
\end{align}
If we define a discretization operation of the form
\begin{equation}
	R^{ij}_{\; \; tx} = M^i_{\; t} \ B^j_{\; x} = \delta(t_i - t) \ \delta(x_j - x) \ ,
\end{equation}
with $x_j = \nicefrac{j}{K}$ for $j \in \left\lbrace 0, 1, ..., K-1 \right\rbrace$, it follows from Eqs.\ \eqref{eq:apiwp} and \eqref{eq:apiwp2} that all Fourier modes are independent and follow IWP processes of the form:
\begin{align}
	&\left.P\left(\left(\begin{matrix}
	\tilde{s}^{i k} \\ \dot{\tilde{s}}^{i k}
	\end{matrix}\right) \right| \left(\begin{matrix}
	\tilde{s}^{(i-1) k} \\ \dot{\tilde{s}}^{(i-1) k}
	\end{matrix}\right) \right) \notag\\ &= \mathcal{G}\left(\left(\begin{matrix}
	\tilde{s}^{i k} \\ \dot{\tilde{s}}^{i k}
	\end{matrix}\right) - \left(\begin{matrix}
	1 & \Delta_i \\
	0 & 1
	\end{matrix}\right) \left(\begin{matrix}
	\tilde{s}^{(i-1) k} \\ \dot{\tilde{s}}^{(i-1) k}
	\end{matrix}\right), \left|\sigma^k\right|^2 \left(\begin{matrix}
	\nicefrac{\Delta_i^3}{3} & \nicefrac{\Delta_i^2}{2} \\
	\nicefrac{\Delta_i^2}{2} & \Delta_i
	\end{matrix}\right)\right),
\end{align}
with $\Delta_i = t_i - t_{i-1}$ and $\tilde{s}^{i k} = (M \tilde{s})^{i k}$. 

As $x_j$ is sampled on a regular grid, from Eq.\ \eqref{eq:fouriersource} we get that
\begin{align}
	\left(\tilde{\bar{s}}^{(c)}\right)^{i k} &= \sum_{n = -\infty}^\infty \left(2 \pi i \left(k+n K\right)\right)^c \tilde{s}^{i (k + n K)} \ , \\
	\left(\dot{\tilde{\bar{s}}}^{(c)}\right)^{i k} &= \sum_{n = -\infty}^\infty \left(2 \pi i \left(k+n K\right)\right)^c \dot{\tilde{s}}^{i (k + n K)}  \ ,
\end{align}
with $k \in \left[-\nicefrac{K}{2}+1, \nicefrac{K}{2}\right]$.

{\bf Proposition:} The random vectors $\tilde{\bar{\mathbf{s}}} = \left(\tilde{\bar{s}}^{(0)}, \tilde{\bar{s}}^{(1)}, ...\right)$ and $\dot{\tilde{\bar{\mathbf{s}}}} = \left(\dot{\tilde{\bar{s}}}^{(0)}, \dot{\tilde{\bar{s}}}^{(1)}, ...\right)$ are Gaussian distributed according to Eq.\ \eqref{eq:pdemarkovprior}.

As the involved discretization operation is a linear operation, it is sufficient to show that the mean and covariance take the proposed form, since $\tilde{s}$ and $\dot{\tilde{s}}$ are itself Gaussian distributed.
For the mean we get that
\begin{align}
	&\left< \left(\tilde{\bar{s}}^{(c)}\right)^{i k} \right> = \sum_{n = -\infty}^\infty \left(2 \pi i \left(k+n K\right)\right)^c \left< \tilde{s}^{i (k + n K)} \right>  \notag \\ &=
	\sum_{n = -\infty}^\infty \left(2 \pi i \left(k+n K\right)\right)^c \left(\tilde{s}^{(i-1) (k + n K)} + \Delta_i \dot{\tilde{s}}^{(i-1) (k + n K)}\right) 
	\notag\\  &= \underbrace{\sum_{n = \-\infty}^\infty \left(2 \pi i \left(k+n K\right)\right)^c \tilde{s}^{(i-1) (k + n K)}}_{= \left(\tilde{\bar{s}}^{(c)}\right)^{(i-1) k}}
	\notag\\ &+ \Delta_i \underbrace{\sum_{n = -\infty}^\infty \left(2 \pi i \left(k+n K\right)\right)^c \dot{\tilde{s}}^{(i-1) (k + n K)} }_{= \left(\dot{\tilde{\bar{s}}}^{(c)}\right)^{(i-1) k}} \notag\\
	&= \left(\tilde{\bar{s}}^{(c)}\right)^{(i-1) k} + \Delta_i \left(\dot{\tilde{\bar{s}}}^{(c)}\right)^{(i-1) k} \ ,
\end{align}
and similarly
\begin{align}
	\left< \left(\dot{\tilde{\bar{s}}}^{(c)}\right)^{i k} \right> &= \sum_{n = -\infty}^\infty \left(2 \pi i \left(k+n K\right)\right)^c \left< \dot{\tilde{s}}^{i (k + n K)} \right> \notag\\ &= \sum_{n = -\infty}^\infty \left(2 \pi i \left(k+n K\right)\right)^c  \dot{\tilde{s}}^{(i-1) (k + n K)} \notag\\
	&= \left(\dot{\tilde{\bar{s}}}^{(c)}\right)^{(i-1) k} \ .
\end{align}

For the equal time covariance we get
\begin{align}
	&\left< \left(\left(\dot{\tilde{\bar{s}}}^{(c)}\right)^{i k} - \left< \left(\dot{\tilde{\bar{s}}}^{(c)}\right)^{i k} \right>\right) \left(\left(\dot{\tilde{\bar{s}}}^{(d)}\right)^{i q} - \left< \left(\dot{\tilde{\bar{s}}}^{(d)}\right)^{i q} \right>\right)^* \right> \notag\\ &= \sum_{n, m = -\infty}^\infty \left(2 \pi i \left(k+n K\right)\right)^c \left(-2 \pi i \left(q+m K\right)\right)^d \times \notag\\ &\underbrace{\left< \left(\dot{\tilde{s}}^{i (k + n K)} - \left<\dot{\tilde{s}}^{i (k + n K)}\right>\right) \left(\dot{\tilde{s}}^{i (q + m K)} - \left<\dot{\tilde{s}}^{i (q + m K)}\right>\right)^* \right>}_{\delta_{n m} \delta_{k q} \left|\sigma^{k+n K}\right|^2 \Delta_i} \notag\\  &=
	\delta_{k q} \ \Delta_i \ (-1)^d \sum_{n = -\infty}^\infty \left(2 \pi i \left(k+n K\right)\right)^{c+d} \left|\sigma^{k+n K}\right|^2  \notag\\ &=
	\delta_{k q} \ \Delta_i \left(\mathbf{D}^k\right)^{c d} \ ,
\end{align}
where we recover the definition of $\mathbf{D}^k$ (Eq.\  \eqref{eq:cov}). An analogous computation of the covariance of $\tilde{\bar{s}}^{(c)}$ yields the same result with $\Delta_i$ being replaced by $\nicefrac{\Delta_i^3}{3}$. Similarly the cross correlation between $\tilde{\bar{s}}^{(c)}$ and its time derivative also results in the same covariance with a pre-factor of $\nicefrac{\Delta_i^2}{2}$.

\end{document}